\documentclass[12pt]{article}
\title{An improved upper bound on the length of the longest cycle of a supercritical random graph}
\author{
Graeme Kemkes\thanks{Partially supported by NSERC CGS and PDF awards. Some of this research appeared in this author's PhD thesis at the University of Waterloo.}\\
{\small Department of Mathematics }\\
{\small University of California, San Diego }\\
{\small La Jolla, California, USA}\\
{\small {\tt gkemkes@math.ucsd.edu }}
\and
Nicholas Wormald\thanks{Supported by the  Canada Research Chairs Program and NSERC.}\\
{\small Department of Combinatorics and Optimization} \\
{\small University of Waterloo}\\
{\small Waterloo ON, Canada}\\
{\small {\tt nwormald@uwaterloo.ca }}
}
\date{}
\usepackage{graphics,graphicx}
  \usepackage{amsmath}

\newcommand{\remove}[1]{}

\def\fullpage {
\addtolength{\topmargin}{-1.5 cm}
\addtolength{\oddsidemargin}{-1.5 cm}
\addtolength{\textwidth}{+3 cm}
\addtolength{\textheight}{+3 cm} }

 \fullpage

       \newcommand{\lab}[1]{\label{#1}}

\catcode`@=11
\@addtoreset{equation}{section}

\catcode`@=12

\def\blackslug{\hbox{\kern1pt\vrule height6pt width4pt  depth1pt\kern1pt}}
\def\qed{\penalty 500\hbox{\quad\blackslug}\ifmmode\else\par
                                            \vskip4.5pt plus3pt minus2pt\fi}
\newtheorem{thm}{Theorem}

\newtheorem{lemma}[thm]{Lemma}
\newtheorem{cor}[thm]{Corollary}
\newtheorem{prop}[thm]{Proposition}

\def\proof{\par\noindent{\bf Proof.\enspace}\rm}

\def\no{\noindent}

\def\commonExp{E^*}
\def\Var{{\bf V}}

\def\eps{\epsilon}
  \def\G{{\cal G}}
  
  \def\P{{\cal P}}

\def\Hc{{\cal H}}

\def\ex{{\bf E}}
\def\pr{{\bf P}}
\def\vd{{\bf d}}
\def\vb{{\bf b}}

\def\xv{{\bf x}}
\def\Yv{{\bf Y}}

\def\Gn5{\G_{n,5}}

\def\Pn5{\P_{n,5}}

\def\Gnm{{\cal G}(n,M)}
\def\multi{{\bf Multi}}
\def\ETb{E_T}

\def\be{\begin{equation}}
\def\ee{\end{equation}}
\def\bea{\begin{eqnarray}}
\def\eea{\end{eqnarray}}
\def\bean{\begin{eqnarray*}}
\def\eean{\end{eqnarray*}}

\newcommand{\eqn}[1]{(\ref{#1})}
\newcommand{\bel}[1]{\be\lab{#1}}

\begin{document}
\maketitle
\begin{abstract}
We improve {\L}uczak's upper bounds on the length of the longest cycle in the random graph $\G(n,M)$ in the ``supercritical phase" where $M=n/2+s$ and
$s=o(n)$ but $n^{2/3}=o(s)$.
The new upper bound is $(6.958+o(1))s^2/n$ with probability $1-o(1)$ as
$n\to\infty$.
Letting $c=1+2s/n$,
the equivalence between $\G(n,p)$ and $\G(n,M)$ implies the same result
for $\G(n,p)$ where $p=c/n$, $c\to 1$, $c-1 =\omega(n^{-1/3})$.
\end{abstract}

\section{Introduction}
The probability space $\Gnm$ of all $n$-vertex graphs with $M$ edges
under the uniform distribution is also known as the uniform random graph model.
It is one of the earliest models of random graphs, originating in a simple
model introduced by Erd\H os~\cite{Erd47}.
We say that $\Gnm$ has a property \emph{asymptotically almost surely}
(abbreviated a.a.s.) if the probability of this event
is $1-o(1)$ as $n\to\infty$.
Much of the interest in this model comes from the study of
its its asymptotically almost sure (also abbreviated a.a.s.) properties
as the dependence of $M$ upon $n$ is varied. This change from a sparse
graph to a dense graph, as $M$ increases more quickly with $n$, is
called the evolution of the random graph.
One important property is
the number $L$ of vertices in the largest component of $\Gnm$. (If there
is more than one component with the maximum number of vertices,
we use the lexicographically first
among largest components.)
When $M=cn/2$ for constant $c$, Erd\H os and R\'enyi~\cite{ER}
showed that the number of vertices in the largest component of $\Gnm$
is a.a.s.\ $O(\log n)$, $\Theta(n^{2/3})$, or $\Theta(n)$ according to whether
$c<1$, $c=1$, or $c>1$, respectively.

Because of this dramatic change in the structure of $\Gnm$, we often call
$M=n/2$ a ``phase transition''. Further research showed that the phase
transition extends throughout the period $M=n/2+cn^{2/3}$ for constant $c$
in the sense that, for this range of $M$,
$L=c'n^{2/3}$ with a distribution over the constant $c'$. As a result,
this range of
$M$ is known as the \emph{critical period}.
For $s=s(n)$ satisfying $n^{2/3}=o(s)$ but $s=o(n)$, the range $M=n/2-s$ is
known as the \emph{subcritical phase} while the range $M=n/2+s$ is known as
the \emph{supercritical phase}.
For $M$ in the supercritical phase, $\Gnm$ a.a.s.\ has a unique
largest component on $(4+o(1))s$ vertices
and every other component has fewer than
$n^{2/3}$ vertices. A ``giant component'' has emerged.

Another well-studied graph property is its \emph{circumference}, the
length of its longest cycle. The circumference $l$ of $\Gnm$ also changes
dramatically during the phase transition, but the way it changes is not entirely understood.
Let $\omega=\omega(n)\to\infty$. When $M=cn/2$ for fixed $c<1$, the
circumference of $\Gnm$ is a.a.s.\ at most $\omega$
(\cite{Bo01}, Corollary~5.8).
In the subcritical phase, the circumference
$l$ a.a.s.\ satisfies $l/\omega < n/s < l\omega$
(\cite{JLR}, Section~5.4). During the critical period $M=n/2+O(n^{2/3})$
it a.a.s.\ satisfies $l/\omega < n^{1/3} < l\omega$
(\cite{JLR}, Section~5.5). But for larger $M$
there are not such good estimates for the circumference. (Of course, when
$M=n(\log n + \log\log n + \omega)/2$ the circumference is a.a.s.\ equal
to $n$ as the graph is a.a.s.\ Hamiltonian~\cite{KoSz83}.)
When $M=cn/2$ for fixed $c>1$, there are several known a.a.s.\ lower bounds
on the circumference of the form $(f(c)+o(1))n$ \cite{FdlV,BFF,F86}. One of the earliest and most significant breakthroughs
was given by Ajtai, Koml\'os and Szemer\'edi \cite{AKS}, who also showed
an equivalence between the problems of finding
paths of length $(f(c)+o(1))n$ and finding
cycles of length $(f(c)+o(1))n$.
{\L}uczak~\cite{Lcyc} has shown that the circumference of $\Gnm$ is
a.a.s.\ between $(16/3+o(1))s^2/n$ and $(7.496+o(1))s^2/n$
when $M=n/2+s$ for $s=o(n)$ and $n^{2/3}=o(s)$.
Kim and the second author~\cite{KW} have improved {\L}uczak's lower bound to
$(6+o(1))s^2/n$.

In this paper we improve upon {\L}uczak's upper bound as follows.
\begin{thm}
\lab{tMain}
Let $M=n/2+s$ with $n^{2/3}=o(s)$ and $s=o(n)$. The circumference of
$\Gnm$ is a.a.s.\ at most $(6.958+o(1))s^2/n$.
\end{thm}

In proving his result,
{\L}uczak
focused on the \emph{core}  and \emph{kernel} of $\Gnm$.
The core of a graph is its maximal subgraph
of minimum degree at least 2.
The \emph{prekernel} of a graph is obtained from the core by throwing away any
cycle components. The \emph{kernel} of a graph is obtained from the prekernel by
replacing each maximal path of degree-2 vertices by a single edge.
We say that a graph is a prekernel (respectively, a kernel) if it is
the prekernel (respectively, kernel) of some graph.

{\L}uczak's main insight was
that, for this range of $M$,
the kernel is much like a random 3-regular graph,
and the core is much like the graph formed from the kernel by randomly subdividing its edges about
$(8+o(1))s^2/n$ times.
A random 3-regular graph a.a.s.\ contains a Hamilton cycle. This gives
a cycle in $\Gnm$ containing about $(2/3)\times (8+o(1))s^2/n = (16/3+o(1))s^2/n$ vertices of the core.
This
was {\L}uczak's lower bound on the circumference.
His upper bound
came from viewing the core as constructed from the kernel
together with a sequence of numbers, summing to $(8+o(1))s^2/n$, describing
how many degree-2 vertices belong on each edge of the kernel. From
probability theory, the sum of the largest
two-thirds of the terms of such a random sequence is at most
$(7.496+o(1))s^2/n$.

We obtain our result by a different, more detailed
study of how a cycle can pass through such a structure.
Our main tool is
the kernel configuration model, introduced in~\cite{PW}  to facilitate
arguments like {\L}uczak's.

Following {\L}uczak's example, it is helpful to put weights on the edges
of the kernel; the weight of an edge tells us how many times the edge should
be subdivided to recover the core. These weights form a random sequence
whose asymptotic properties we investigate in
Section~2. In particular, we show that any bounded number of terms in such
a sequence behave like independent random variables with exponential
distribution. We also show that when a function of a bounded number of these terms is
summed over many sets of such terms, the result is concentrated about its
expected value.
 These properties
are needed in Section~3 where we establish an a.a.s.\ upper bound
on the weight of the heaviest cycle in a pseudograph
 with random edge weights. The upper bound
is expressed in terms of
a family of constants, some of which we explicitly calculate in Section~4.
In Section~5 we prove an a.a.s.\ upper bound on the circumference of a
random prekernel
with a degree sequence that resembles a random 3-regular graph with
subdivided edges.
In Section~7 we use this result to prove Theorem~1
after, in Section~6, establishing that the
degree sequence of the prekernel of $\Gnm$ indeed shows the required resemblance.
\section{Random sequences}
Let $\Omega$ be the probability space, equipped with the uniform distribution,
of all sequences of $m$
positive integers $(X_1, X_2, \ldots, X_m)$ summing to $N$.
We are interested in the asymptotic value of certain functions of these
random variables. Letting $\omega = \omega(N)\to\infty$, our asymptotics are
in terms of $N\to\infty$, uniformly
over all $m$ satisfying $\omega < m < N/\omega$.
For the rest of the paper
we write $\mu = N/m$.

Our first result tells us the expected value of certain functions of
$X_1, X_2, \ldots, X_j$ for $j$ bounded.

\begin{lemma}
\lab{lExp}
Let $g$ be a
nonnegative
integrable function of a bounded number $j$ of
nonnegative variables.
Suppose that for some $C$ and $d$,
$g( x_1, \ldots, x_j ) \le C (x_1 + \cdots + x_j)^d$ for all
$x_1,  \ldots, x_j$.
Then,
\[
\ex\left[g\left(\frac{X_1}{\mu},\ldots,\frac{X_j}{\mu}\right)\right]
=
\int_0^\infty \cdots \int_0^\infty g(x_1, \ldots, x_j)
e^{ -x_1 - x_2 - \cdots -x_j } dx_1 \cdots dx_j
+ o(1).
\]
\end{lemma}

Since the $X_i$ are identically distributed, the above theorem
also holds when
$( X_1, \ldots, X_j )$
is replaced by $(X_{\sigma(1)}, \ldots, X_{\sigma(j)})$
for any $j$ distinct $\sigma(1), \ldots, \sigma(j)$ in
$\{1, 2, \ldots, m\}$. Furthermore, the error represented by $o(1)$
is independent of $\sigma$.

The next result states that when such a function
is summed over
$\sigma$ in a sufficiently rich family,
the sum is asymptotically almost surely (a.a.s.) concentrated about
its expected value.

\begin{lemma}
\lab{lConc}
Let $f$ be a nonnegative integrable function of a bounded number $k$ of
nonnegative variables.
Suppose that for some $C$ and $d$,
$f( x_1, \ldots, x_k ) \le C (x_1 + \cdots + x_k)^d$ for all
$x_1,  \ldots, x_k$.
Define the constant
\[
\commonExp :=
\int_0^\infty \cdots \int_0^\infty f(x_1, \ldots, x_k)
e^{ -x_1 - x_2 - \cdots -x_k } dx_1 \cdots dx_k.
\]
and assume $\commonExp > 0$.
Let $S$ be a set of $k$-tuples with entries from
$\{1, 2, \ldots, m\}$,
with each $k$-tuple having distinct components.
Let $I=I(S) \in S \times S$ be the pairs of tuples which intersect; that is,
\[
I = \{ (\sigma, \tau) \in S \times S \mid
\{\sigma(1), \ldots, \sigma(k)\} \cap \{\tau(1), \ldots, \tau(k)\} \ne \emptyset \}.
\]
If $|I|=o(|S|^2)$ then
\[
\sum_{\sigma \in S}
  f\left( \frac{X_{\sigma(1)}}{\mu}, \ldots, \frac{X_{\sigma(k)}}{\mu} \right)
=
(\commonExp + o(1))
|S|
\]
a.a.s.; that is, with probability $1-o(1)$. Furthermore, the
$o(1)$ terms may be bounded independently of $S$.
\end{lemma}
These types of concentration results are often proved using martingales
or inequalities like Talagrand's; however, because we are aiming for such
a coarse result, a simple application of Chebyshev's inequality will suffice
for the proof.



\subsection{Distribution of terms}
In this section we establish some preliminary results about the distribution
of the positive terms $X_1, X_2, \ldots, X_j$ for bounded $j$.
It is an exercise in basic counting to show that the number of sequences in
$\Omega$ is $\binom{N-1}{m-1}$. It immediately follows that for positive integers $t_1, t_2, \ldots,
t_j$, the number of sequences in $\Omega$ with $X_1=t_1$, $X_2=t_2$, $\ldots$,
$X_j=t_j$ is
\[
B(t) := \binom{N-1-t}{m-j-1}
\]
where $t = t_1 + t_2 + \cdots + t_j$.

\begin{prop}
\lab{propB}
Let $x$ satisfy
$x<\sqrt{m}/\omega$ and $x<\mu/\omega$.
For positive integers $t \le x\mu$ we have
\[
\frac{B(t)}{|\Omega|} = (1+O(\omega^{-1})) \mu^{-j} e^{-t/\mu}.
\]
\end{prop}
\proof
\begin{eqnarray*}
\frac{B(t)}{|\Omega|}
&=& \binom{N-1-t}{m-j-1} \binom{N-1}{m-1}^{-1} \\
&=& \bigg(\prod_{i=1}^{m-1}\frac{ N-t-i}{N-i} \bigg) \bigg( \prod_{i=1}^j\frac{m-i}{N-t-m+i}\bigg) \\
&=& \bigg(\prod_{i=1}^{m-1}\bigg(1-\frac{t}{ N-i}\bigg)\bigg)
\mu^{-j}\big(1+O(m^{-1})+O(t/N)+O(m/N)\big).
\eean
Since $m>\omega$, $x<m/\omega$  and $m<N/\omega$, the error term is $O(1/\omega)$. Also
\bean
1-\frac{t}{ N-i}
&=&  1-\frac{t}{ N}\bigg(1+O\bigg(\frac{i}{N}\bigg)\bigg)\\
&=&  1-\frac{t}{ N} (1+O (\mu^{-1}));\\
\prod_{i=1}^{m-1}\bigg(1-\frac{t}{ N-i}\bigg)
&=&  e^{-t/\mu}\big(1+O(tm^2/N^2+mt^2/N^2)\big)\\
&=&   e^{-t/\mu}
      \big(1-O(\omega^{-1})\big).  \qed
\end{eqnarray*}

\begin{cor}
\lab{corDistMain}
Let $x>0$ be fixed. For any positive integers $t_1, t_2, \ldots, t_j$ summing to
$t \le x\mu$ we have
\[
\pr[X_1=t_1, X_2=t_2, \ldots, X_j=t_j] = (1+O(\omega^{-1})) \mu^{-j} e^{-t/\mu}.
\]
\end{cor}

Next we bound the probability of larger terms.

\begin{lemma}
\lab{lemmaDistTail}
Let $x>0$ be fixed. For positive integers $t_1, t_2, \ldots, t_j$ summing to
$t \ge x\mu$ we have
\[
\pr[X_1=t_1, X_2=t_2, \ldots, X_j=t_j]
< 2 \mu^{-j} e^{-x} \left( 1 - \frac{1}{2\mu} \right)^{t-x\mu}
\]
when $N$ is sufficiently large.
\end{lemma}
\proof
If $B(t)=0$ then the required probability is zero and we are done.
Otherwise, $B(i)$ is nonzero for all
positive integers $i\le t$ and the probability which we must estimate is
\begin{eqnarray*}
\frac{B(t)}{|\Omega|}
&=& |\Omega|^{-1} B(\lfloor x\mu \rfloor)
      \prod_{i=\lfloor x\mu \rfloor + 1}^t \frac{B(i)}{B(i-1)}.
\end{eqnarray*}
By Proposition~\ref{propB}, the product of the first two terms is
$(1+O(\omega^{-1}))\mu^{-j} e^{-\lfloor x\mu \rfloor / \mu}$.
This is less than $2 \mu^{-j} e^{-x}$ when $N$ is sufficiently large.
To bound the remaining product, we estimate the ratio
\begin{eqnarray*}
\frac{B(i)}{B(i-1)}
&=& \frac{\binom{N-1-i}{m-j-1}}{\binom{N-1-i+1}{m-j-1}} \\
&=& 1 - \frac{m-j-1}{N-i} \\
&<& 1 - \frac{m-j-1}{N} \\
&<& 1 - \frac{m/2}{N}
\end{eqnarray*}
where the last inequality holds for $N$ sufficiently large. So, for $N$
sufficiently large,
$$
      \prod_{i=\lfloor x\mu \rfloor + 1}^t \frac{B(i)}{B(i-1)}
\ <\  \left(1 - \frac{1}{2\mu} \right)^{t-\lfloor x\mu \rfloor}
\ \le\ \left(1 - \frac{1}{2\mu} \right)^{t- x\mu}.
$$
The result follows.
\qed

\subsection{Proof of Lemma~\ref{lExp}}
By the definition of expected value, we have
\[
\ex\left[g\left(\frac{X_1}{\mu},\ldots,\frac{X_j}{\mu}\right)\right]
=
\sum g\left(\frac{t_1}{\mu},\ldots,\frac{t_j}{\mu}\right)
   \pr[X_1=t_1, X_2=t_2, \ldots, X_j=t_j]
\]
where the sum is over all positive integer $j$-tuples
$t_1, t_2, \ldots, t_j$.

Fix $x>0$. Let us split the sum into two parts, $S_1(x)$ being the sum over
$j$-tuples where each $t_i < x\mu$, and $S_2(x)$ being the remainder.
We will show that, as $N \to \infty$,
\[
S_1(x) \to
\int_0^x \cdots \int_0^x g(x_1, \ldots, x_j)
e^{ -x_1 - x_2 - \cdots -x_j } dx_1 \cdots dx_j,
\]
while
\[
|S_2(x)| < K e^{-x/2}
\]
for some constant $K$. As $x$ grows,
$|S_2(x)|$ approaches 0
and $S_1(x)$ is
nonnegative and nondecreasing
since $g$ is nonnegative.
So, taking $x \to \infty$ proves the lemma.

We begin by estimating $S_1(x)$. These terms have each $t_i \le x\mu$, so
we use Corollary~\ref{corDistMain} to estimate the probabilities as follows.
\begin{eqnarray*}
S_1(x)
&=& \sum_{t_1<x\mu} \cdots \sum_{t_j<x\mu}
 g\left(\frac{t_1}{\mu},\ldots,\frac{t_j}{\mu}\right)
   \pr[X_1=t_1, X_2=t_2, \ldots, X_j=t_j]  \\
&=& \sum_{t_1<x\mu} \cdots \sum_{t_j<x\mu}
 g\left(\frac{t_1}{\mu},\ldots,\frac{t_j}{\mu}\right)
 (1+O(\omega^{-1})) \mu^{-j} e^{-(t_1+\cdots+t_j)/\mu}.
\end{eqnarray*}
Since $O(\omega^{-1})$ is independent of the $t_i$, this becomes
\[
 (1+O(\omega^{-1}))\sum_{t_1<x\mu} \cdots \sum_{t_j<x\mu}
 g\left(\frac{t_1}{\mu},\ldots,\frac{t_j}{\mu}\right)
  \mu^{-j} e^{-(t_1+\cdots+t_j)/\mu}.
\]
Letting $M=x\mu$ we get
\[
 (1+O(\omega^{-1}))\sum_{t_1<M} \cdots \sum_{t_j<M}
 g\left(t_1\frac{x}{M},\ldots,t_j\frac{x}{M}\right)
  e^{-(t_1+\cdots+t_j)x/M}
  \left(\frac{x}{M}\right)^{j}.
\]
As $N\to\infty$ we have $M\to\infty$ and this expression becomes the
 Riemann integral
\[
\int_0^x \cdots \int_0^x g(x_1, \ldots, x_j)
e^{ -x_1 - x_2 - \cdots -x_j } dx_1 \cdots dx_j
\]
as required.

The terms of the sum $S_2(x)$ are indexed by $j$-tuples $t_1, t_2, \ldots, t_j$
with at least one $t_i \ge x\mu$. Consider such a term, and let
$t = t_1+t_2+\cdots+t_j$. For $N$ sufficiently large, the absolute value
of the term is
\[
 g\left(\frac{t_1}{\mu},\ldots,\frac{t_j}{\mu}\right)
   \pr[X_1=t_1, X_2=t_2, \ldots, X_j=t_j]
<
C\left(\frac{t}{\mu}\right)^d
 2 \mu^{-j} e^{-x} \left( 1 - \frac{1}{2\mu} \right)^{t-x\mu}
\]
by the hypotheses about $g$ and Lemma~\ref{lemmaDistTail}.
The number of terms in $S_2(x)$ indexed by $j$-tuples summing to $t$ is at
most $\binom{t-1}{j-1} \le (t+j)^{j-1}  \le (2t)^{j-1}$ for $N$ (and hence
$t$) sufficiently large. Thus, for $N$ large, we have
\begin{eqnarray*}
|S_2(x)|
&<& \sum_{t\ge x\mu} (2t)^{j-1}
C\left(\frac{t}{\mu}\right)^d
 2 \mu^{-j} e^{-x} \left( 1 - \frac{1}{2\mu} \right)^{t-x\mu} \\
&=&
2e^{-x}\frac{C2^{j-1}}{\mu^{j+d}}
  \left( 1 - \frac{1}{2\mu} \right)^{-x\mu}
\sum_{t\ge x\mu} t^{j+d-1}
  \left( 1 - \frac{1}{2\mu} \right)^{t}.
\end{eqnarray*}
The factor $( 1 - 1/(2\mu))^{-x\mu}$ approaches $e^{x/2}$ as $N\to\infty$.
The remaining sum is
\begin{eqnarray*}
\sum_{t\ge x\mu} t^{j+d-1}
  \left( 1 - \frac{1}{2\mu} \right)^{t}
&\le&
\sum_{t\ge 0} (t+1)(t+2)\cdots (t+j+d-1)\left(1-\frac{1}{2\mu}\right)^t\\
&=&
(j+d-1)!(2\mu)^{j+d}
\end{eqnarray*}
using the Maclaurin series expansion
$k!(1-x)^{-k-1}=\sum_{t\ge 0}(t+1)(t+2)\cdots (t+k)x^t$.
Combining this with the previous results, we get the desired estimate.
This proves the lemma.
\qed

\subsection{Proof of Lemma~\ref{lConc}}
For each $\sigma$ in $S$, define the random variable
$Y_{\sigma} := f( X_{\sigma(1)}/\mu, \ldots, X_{\sigma(k)}/\mu )$.
As we remarked after Lemma~\ref{lExp}, each of these variables has
the same distribution as the random variable
$Y_{1} := f( X_{1}/\mu, \ldots, X_{k}/\mu )$. In particular, the
expected value is the constant $\commonExp$, up to an additive error of $o(1)$.
We will establish the concentration of the random variable
$Z := \sum_{\sigma \in S} Y_{\sigma}$
by showing that the variance $\Var[Z]$ is $o((\ex Z)^2)$. The lemma then
follows by Chebyshev's inequality.

We begin by estimating
\begin{eqnarray*}
(\ex Z)^2 &=& \sum_{(\sigma, \tau) \in S \times S}
          \ex Y_\sigma \ex Y_\tau  \\
&=& \sum_{(\sigma, \tau) \in S \times S}
          (\commonExp + o(1))
          (\commonExp + o(1)) \\
&=& \sum_{(\sigma, \tau) \in S \times S}
          \Theta(1) \\
&=& \Theta(|S|^2) \\
\end{eqnarray*}
 (using the lower bound assumed on $\commonExp$ in the lemma).
We can write the variance as
\begin{eqnarray*}
\Var[Z] &=& \ex[Z^2] - (\ex Z)^2 \\
&=& \sum_{(\sigma, \tau) \in S \times S} ( \ex[Y_\sigma Y_\tau] - \ex Y_\sigma \ex Y_\tau ) \\
&=& \sum_{(\sigma, \tau) \in I}
         ( \ex[Y_\sigma Y_\tau] - \ex Y_\sigma \ex Y_\tau )
+ \sum_{(\sigma, \tau) \in (S \times S) \setminus I }
         ( \ex[Y_\sigma Y_\tau] - \ex Y_\sigma \ex Y_\tau ).
\end{eqnarray*}
To study the terms of the second sum, let
$(\sigma, \tau)\in (S\times S)\setminus I$. By Lemma~\ref{lExp},
we have
\begin{eqnarray*}
&&\ex[ Y_\sigma Y_\tau ] \\
&=&
 \ex\left[
    f\left( \frac{X_{\sigma(1)}}{\mu}, \ldots, \frac{X_{\sigma(k)}}{\mu}\right)
    f\left( \frac{X_{\tau(1)}}{\mu}, \ldots, \frac{X_{\tau(k)}}{\mu}\right)
\right]\\
&=&
 \int_0^\infty \cdots \int_0^\infty
    f( x_1, \ldots, x_k )
    f( x_{k+1}, \ldots, x_{2k} )
    e^{ -x_1 -  \cdots -x_{2k} } dx_1 \cdots dx_{2k} + o(1)\\
&=&
\left( \int_0^\infty \cdots \int_0^\infty
    f( x_1, \ldots, x_k )
    e^{ -x_1 - x_2 - \cdots -x_{k} } dx_1 \cdots dx_{k} \right)^2 + o(1)\\
&=&
\ex Y_\sigma \ex Y_\tau + o(1) \\
\end{eqnarray*}
where $o(1)$ is independent of $\sigma$ and $\tau$. So the second sum is
$o(|S|^2)$.
To study the terms of the first sum, we can be more crude.
By Lemma~\ref{lExp} and the remark following it, we know  that each
$\ex[ Y_\sigma Y_\tau ]$ and $\ex Y_\sigma \ex Y_\tau$ depends only on
the tuple positions where $\sigma$ and $\tau$ intersect, and each value
is $O(1)$. So the first sum is $O(|I|)$, which is $o(|S|^2)$ by hypothesis.
Combining the two sums, we see that
the variance of $Z$ is $o(|S|^2)$, which is $o((\ex Z)^2)$, as
required.
\qed

\section{Heavy cycles in a weighted pseudograph}
In the introduction we saw that the problem of bounding the circumference
of $\Gnm$ is connected to the problem of bounding the weight of the
heaviest cycle in a certain edge-weighted graph.
In this section we study a graph, technically a pseudograph since it may have loops and/or multiple edges, whose $m$ edges are randomly
weighted by positive
integers summing to $N$. The sequence of weights is chosen uniformly at
random from among all such sequences. Equivalently, we can think of the
weights as being generated by the following random process applied to make
a sequence of pseudographs, beginning with the given one. At each step,
choose an edge
uniformly at random from the current pseudograph and subdivide the edge into
two edges.
Repeat the procedure until the resulting pseudograph has exactly $N$
edges.  For each edge in the original pseudograph, define its
weight to be the number of edges into which it has been subdivided. These
weights form a sequence of $m$ positive integers summing to $N$. There
are exactly $(N-m)!$ ways that the process can form a given sequence, so the
sequence is chosen
uniformly at random from among all such sequences.
Another random process for generating the weights initially gives a weight
of 1 to each edge, then selects an edge at random with probability
proportional to the weight of the edge and increments the weight of the
selected edge by 1. The
selection and incrementing is repeated until the total weight is $N$. It is
easy to see that this process is equivalent to the previous one.

Given a subgraph of an edge-weighted pseudograph, we define the weight
of the  subgraph to be the sum of the weights on its edges.
To establish an upper bound for the weight of a cycle in a large
pseudograph, we will consider
the intersection of the cycle with small trees in the pseudograph. The
intersection of the cycle and the small tree will form a set of
vertex-disjoint paths
which begin and end at leaf vertices of the tree. We will use the maximum-weight
set of such vertex-disjoint paths
to bound the weight of the intersection. This motivates the
following definitions.

Fix an integer $k\ge 2$. A \emph{biased tree} $T$ on $k$ edges is a tree on $k$ edges with each
non-leaf vertex
having degree 3 and each edge $e_i$ having a nonnegative number $b_i$ called its bias.
We may assume that the sum of the biases $\vb=(b_1, b_2, \ldots, b_k)$ is 1.

Let $\P$ be the set of all maximal subgraphs of $T$ which are a union of
 vertex-disjoint paths which begin and end at leaf vertices.
 Define the function
\[
f_T(x_1, x_2, \ldots, x_k)
= \max_{P\in \P} \sum_{i: e_i\in E(P)} b_i x_i
\]
and the constant
\bel{eETb}
\ETb =
\int_0^\infty \int_0^\infty \cdots \int_0^\infty
f_T(x_1, x_2, \ldots, x_k)
e^{-x_1-x_2- \cdots -x_k}
dx_1 dx_2 \cdots dx_k.
\ee
If $x_1, x_2, \ldots, x_k$ are weights on the edges of $T$, we can think
of $f$ as the maximum ``biased weight'' of any graph in $\P$.

We say that the positive constant $c^*$ is \emph{$k$-admissible}
if $\ETb < c^*$ for some biased tree $T$ on $k$ edges.

\begin{lemma}
\lab{lemmaHeavy}
Fix an integer $k\ge 2$. Let the positive number $c^*$ be $k$-admissible.
Let $G=G(n)$ be a pseudograph on $v=v(n)\to\infty$ (as $n\to\infty$)
vertices and $m=m(n)$ edges with minimum degree at least 3.
Suppose the subgraph $B$ of $G$ induced by
cycles of length at most $k$ (including loops and parallel edges)
and edges incident to
vertices of degree greater than 3 satisfies $|E(B)| = o(v)$.
Let $N=N(n)$ be a positive integer satisfying $m=o(N)$.
On the edges of $G$ put weights, a sequence chosen uniformly at random
from among all sequences of $m$ positive integers summing to
$N$.  Then, the heaviest cycle in $G$ has weight a.a.s.\ at most
$c^* N$.
\end{lemma}
\proof
Denote the edges of $G$ by $w_1, w_2, \ldots, w_m$ and their
random weights by $X_1, X_2, \ldots, X_m$.  We estimate $m$ by recalling
that in any graph the sum of the vertex degrees equals twice the number
of edges.
Since $G$ has minimum degree at least 3, we have $2m\ge 3v$. Since $G$ has only $o(v)$
edges incident to vertices of degree greater than 3, we have $2m\le 3v+o(v)$.
Thus $m\sim 3v/2$.

For a subgraph $S$ of $G$, define its \emph{$k$-neighbourhood} $\Gamma(S)$ to be
the subgraph of $G$ reachable from $S$ by paths of length at most
$k$. Recalling that the subgraph $B$ of $G$ contains all edges incident with
vertices of degree greater than 3, its $k$-neighbourhood satisfies
$|E(\Gamma(B))| \le 2^k |E(B)| = o(v)$.

Let $C$ be a cycle in $G$. Its weight $wt(C)$ is
\begin{eqnarray*}
wt(C) &:=& \sum_{j=1}^m X_j I(w_j\in E(C))\\
     &=&  \sum_{j:w_j\in E(\Gamma(B))} X_j I(w_j\in E(C))
        + \sum_{j:w_j\not\in E(\Gamma(B))} X_j I(w_j\in E(C))
\end{eqnarray*}
where $I(\alpha)$ is the indicator function equal to 1 if $\alpha$
is true and 0 otherwise.
The expected value of each $X_j$ is $\mu=N/m$, so the first sum has expected
value at most $|E(\Gamma(B))|N/m = o(vN/m) =o(N)$ since $m\sim 3v/2$. It follows by Markov's inequality
that the first sum is a.a.s.\ $o(N)$. Thus a.a.s.,
\bel{eqwC}
wt(C) = o(N) + \sum_{j:w_j\not\in E(\Gamma(B))} X_j I(w_j\in E(C)).
\ee

Since $c^*$ is $k$-admissible, there is a biased tree $T$ on $k$ edges
with $\ETb < c^*$.
We will study the copies of $T$ in the graph $G\setminus B$. Let $S$ be the set
of 1-1 homomorphisms $\sigma$ mapping $T$ to $G\setminus B$. Since $k\ge 2$,
each $\sigma$ is uniquely defined by the mapping it induces between the edge
sets.
We write $\sigma = (\sigma(1), \sigma(2), \ldots, \sigma(k))$ and interpret
$\sigma(i)=j$ to mean that $\sigma$ maps edge $e_i$ of $T$ to edge $w_j$ of $G$.

Consider the random variable
\[
Z=\sum_{\sigma\in S}\sum_{i=1}^k b_i X_{\sigma(i)} I(w_{\sigma(i)}\in E(C)).
\]
Expressing $Z$ in terms of the edges of $G$ we may write
\begin{eqnarray*}
Z &=& \sum_{j=1}^m \sum_{\sigma\in S} \sum_{i=1}^k b_i X_j I(w_j\in E(C)) I(\sigma(i)=j)\\
&=& \sum_{j=1}^m \sum_{i=1}^k b_i X_j I(w_j\in E(C))
|\{ \sigma\in S \mid \sigma(i)=j\}|.
\end{eqnarray*}
For each edge $w_j$ of $G$ not in $E(\Gamma(B))$ the $k$-neighbourhood of $w_j$ is the depth-$k$ tree with
internal vertices of degree 3. Thus,
$|\{ \sigma\in S \mid \sigma(i)=j\}|$
equals some constant independent of $j$. In fact, this constant is
a number $a$, independent of $i$, because any $\sigma$ in this set is
determined by choosing one of the 2 ways to embed $e_i$ onto $w_j$ and then,
moving outward from $e_i$, making one binary choice for each non-leaf vertex
of $T$. On the other hand, for an edge $w_j\in E(\Gamma(B))$,
$|\{ \sigma\in S \mid \sigma(i)=j\}|$
is at most $a$ (by the same argument, recalling that some
choices are impossible because $\sigma$ maps into $G\setminus B$), so we have
\begin{eqnarray*}
Z &=& \sum_{j:w_j\in E(\Gamma(B))} O(X_j)
   +\sum_{j:w_j\not\in E(\Gamma(B))} \sum_{i=1}^k a b_i X_j I(w_j\in E(C))\\
&=& \sum_{j:w_j\in E(\Gamma(B))} O(X_j)
   +\sum_{j:w_j\not\in E(\Gamma(B))}  a  X_j I(w_j\in E(C))\\
\end{eqnarray*}
since $\sum_{i=1}^k b_i=1$.
The first sum has $o(v)$ terms, each having expected value $O(N/m)$,
so the sum is a.a.s.\ $o(vN/m)=o(N)$ by Markov's inequality.
 We now have a.a.s.\
\[
Z = o(N) + a \sum_{j\not\in E(\Gamma(B))} X_j I(w_j\in E(C)).
\]
Combining this result with \eqn{eqwC} we get a.a.s.\
\bel{eqwC2}
wt(C) = \frac{1}{a}Z + o(N).
\ee

Returning to the definition of $Z$, we notice that the inner sum
is the ``biased weight'' of the edges of $C$ passing
through the copy of $T$ given by $\sigma$. These edges must form
vertex-disjoint paths beginning and ending at leaves of the copy of $T$
given by $\sigma$, so this sum is at most
 $f_{T}(X_{\sigma(1)}, X_{\sigma(2)}, \ldots, X_{\sigma(k)})$.  So
\[
Z\le\sum_{\sigma\in S} f_{T}(X_{\sigma(1)}, X_{\sigma(2)}, \ldots, X_{\sigma(k)}).
\]
We will estimate this sum by applying Lemma~\ref{lConc} to
\[
\frac{1}{\mu} \sum_{\sigma\in S} f_{T}(X_{\sigma(1)}, X_{\sigma(2)}, \ldots, X_{\sigma(k)})
=
 \sum_{\sigma\in S} f_{T}\left(\frac{X_{\sigma(1)}}{\mu}, \frac{X_{\sigma(2)}}{\mu}, \ldots, \frac{X_{\sigma(k)}}{\mu}\right).
\]
To verify the hypotheses of Lemma~\ref{lConc} we first note that
$f_T(x_1, x_2, \ldots, x_k)$ is nonnegative, piecewise linear (and hence integrable), and
bounded above by $x_1+x_2+\cdots+x_k$.
We estimate $|S|$ by
\begin{eqnarray*}
|S|
&=& \sum_{j:w_j\in E(G)} |\{\sigma\in S \mid \sigma(1)=j\}| \\
&=& \sum_{j:w_j\in E(\Gamma(B))} |\{\sigma\in S \mid \sigma(1)=j\}|
+ \sum_{j:w_j\not\in E(\Gamma(B))} |\{\sigma\in S \mid \sigma(1)=j\}| \\
&=& o(v) + \sum_{w_j\not\in E(\Gamma(B))} a \\
&=& o(v) + (1+o(1))m a \\
&\sim& ma
\end{eqnarray*}
using $o(v)=o(m)$.
To estimate the cardinality of the set $I$ of pairs
$(\sigma,\tau)\in S\times S$ for which
$\sigma$ and $\tau$ represent intersecting copies of $T$, consider any
edge $f$ in $G\setminus B$. As we have seen previously, there are at most $a$
copies of $T$ using $f$. So, a crude upper bound for $|I|$ is
$a^2|E(G\setminus B)| \le a^2 m$, giving us $|I|=o(|S|^2)$ as required.
Recalling the definition of $\ETb$ from \eqn{eETb}
we may apply Lemma~\ref{lConc} and conclude a.a.s.\
\begin{eqnarray*}
\frac{1}{\mu}Z &\le&  (\ETb+o(1))|S| \\
   &=& (\ETb+o(1))am.
\end{eqnarray*}
Combining this with Equation~\eqn{eqwC2} we get a.a.s.\
\begin{eqnarray*}
wt(C) &\le& \frac{1}{a}\mu (\ETb+o(1))am + o(N) \\
 &=& (\ETb+o(1))N + o(N) \\
 &<& c^* N
\end{eqnarray*}
as required.
\qed
{\bf Remark 1.} A random 3-regular graph a.a.s.\ satisfies all
of the hypotheses of
Lemma~\ref{lemmaHeavy}. The lemma thus gives an upper bound which holds a.a.s.\
on the weight of the heaviest cycle in a randomly-weighted random 3-regular
graph.

{\bf Remark 2.} There are essentially two ingredients in the proof
of Lemma~\ref{lemmaHeavy}. The first ingredient is a method for bounding the
weight of a cycle in a large edge-weighted 3-regular subgraph. The second
ingredient is the argument that the
weight of the heaviest cycle does not change much when the remainder of
the graph is included.
This second ingredient is implicit in {\L}uczak's proof of his upper bound
on the circumference of $\Gnm$ in the supercritical phase~\cite{Lcyc}.
It is the first ingredient that is the new contribution.

{\bf Remark 3.} For the task of bounding the weight of a cycle in
a large edge-weighted 3-regular subgraph, one might suggest investigating
the weight of the least-weight matching. Certainly the
complement of a Hamilton cycle in a 3-regular graph forms a perfect matching.
But, in general, the maximum-weight cycle is not necessarily Hamiltonian.
Thus, its removal from the graph does not always form a perfect matching.

\section{Computing $\ETb$}

Recall the definitions of $f_T(\xv)=f_T(x_1, x_2, \ldots, x_k)$
and $\ETb$ from \eqn{eETb}.
In the previous section we saw that the weight of the heaviest cycle in a
certain edge-weighted pseudograph can be bounded in terms of $\ETb$ for any biased
tree $T$.
In this section we compute
the value of $\ETb$ for a few specific biased trees $T$.
For some trees $T$ we also state the biases $\vb$ which make
$\ETb$ as small as possible.

\begin{prop}
Let $T$ be the biased tree on two degree-3 vertices and four leaf vertices
with bias $b$ on the edges incident to leaves and bias $1-4b$ on the remaining
edge,
where $b$ is the unique zero of $105b^3-90b^2+24b-2$ on $0<b<1/4$.  Then
\[
E_{T} = \frac{4(1-3b)(5b^2-5b+1)}{(7b-2)^2}
\]
which lies in $(0.8797, 0.8798)$ and hence $c^*=0.8798$ is
$k$-admissible for $k=5$.
\end{prop}
It can be shown that, for this tree,
no other choice of biases $\vb$ yields a lower value of $E_{T}$.
See~\cite{K} for details.
\proof
We begin by letting $v$ and $w$ denote the two non-leaf vertices of $T$.
 Let $e_1$ and $e_2$ denote the two edges which are
each incident to $v$ and a leaf. Denote by $e_3$ the edge joining $v$ to $w$,
and christen the other two edges as $e_4$ and $e_5$.
Under this ordering the biases are $\vb = (b,b,1-4b,b,b)$.

To evaluate the integral $E_{T}$ we exploit some of its symmetry. It suffices
to
integrate over only nonnegative $x_1, x_2, x_3, x_4, x_5$ satisfying
$x_1\le x_2$ and $x_4\le x_5$ and multiply the final result by 4.
For such points, only two of the $P\in\P$ can attain the maximum in the
definition of $f_{T}$, giving us
\[
f_T(\xv) =
f_{T}(x_1,x_2,x_3,x_4,x_5) =
\max \{ bx_1+bx_2+bx_4+bx_5, bx_2 + (1-4b)x_3 + bx_5 \}.
\]
We split the region of integration into two parts, according to whether
\[
 bx_1+bx_2+bx_4+bx_5 \ge bx_2 + (1-4b)x_3 + bx_5,
\]
i.e.\ $b(x_1+x_4)/(1-4b) \ge x_3$. The integrals are
\[
\int_{x_2=0}^\infty
\int_{x_5=0}^\infty
\int_{x_1=0}^{x_2}
\int_{x_4=0}^{x_5}
\int_{x_3=0}^{\frac{b(x_1+x_4)}{1-4b}}
b(x_1+x_2+x_4+x_5)
e^{-x_1-x_2-x_3-x_4-x_5}
dx_1 dx_2 \cdots dx_5
\]
and
\[
\int_{x_2=0}^\infty
\int_{x_5=0}^\infty
\int_{x_1=0}^{x_2}
\int_{x_4=0}^{x_5}
\int_{x_3=\frac{b(x_1+x_4)}{1-4b}}^\infty
(bx_2+(1-4b)x_3+bx_5)
e^{-x_1-x_2-x_3-x_4-x_5}
d\xv
\]
which, when evaluated, added together, and multiplied by $4$, give us
\[
E_{T} = \frac{4(1-3b)(5b^2-5b+1)}{(7b-2)^2}.
\]
The result follows by simple computations.
\qed

\begin{prop}
Let $T$ be the biased tree on three degree-3 vertices and five leaf vertices
with bias $b$ on edges incident to leaves and bias $(1-5b)/2$ on
the other edges, where $b$ is the unique zero of
$-3993b^4+2765b^3+1452b^5-804b^2+105b-5$ on $0<b<1/5$.
Then
\[
\ETb = \frac{726b^4-601b^3+245b^2-55b+5}{5(1-b)(4b-1)^2}
\]
which lies in $(0.8741, 0.8742)$ and hence $c^*=0.8742$ is
$k$-admissible for $k=7$.
\end{prop}
For the tree in the above proposition, it can be shown~\cite{K} that
no other choice of biases $\vb$ yields a lower value of $E_{T}$.
\proof
We may view $T$ as the complete binary tree on six edges with
one additional edge $e_1$
joining the root to an additional vertex. Denote by $e_2$ and $e_3$ the other
edges incident to the root. Denote by $e_4$ and $e_5$ the edges incident
with $e_2$. Denote by $e_6$ and $e_7$ the two edges incident with $e_3$.
Under this ordering, the biases are
\[
\vb = \left( b, \frac{1-5b}{2}, \frac{1-5b}{2}, b, b, b, b \right).
\]

By symmetry we may compute $E_{T}$ by
integrating over only
$x_4\le x_5$ and $x_6\le x_7$ and multiplying the final result by 4.
In this range, $f_T$ is the maximum of four expressions,
\begin{enumerate}
\item $b(x_4+x_5+x_6+x_7)$,
\item $bx_1+(1-5b)x_2/2+b(x_5+x_6+x_7)$,
\item $(1-5b)(x_2+x_3)/2+b(x_5+x_7)$, and
\item $bx_1+(1-5b)x_3/2+b(x_4+x_5+x_7)$.
\end{enumerate}

To compute the integral,
the region of integration is divided into four parts, according to which
of the above expressions gives the maximum. We present the details for the
first part only.

The first expression exceeds the other three if and only if
$bx_4 > bx_1 + (1-5b)x_2/2$ and $bx_6 > bx_1 + (1-5b)x_3/2$. To express the integral over
this part as an iterated integral, we divide the part into two regions,
according to whether $x_4>x_6$ or not. The region on which $x_4>x_6$
gives the integral
\[
\int_{x_4=0}^\infty
\int_{x_6=0}^{x_4}
\int_{x_1=0}^{x_6}
\int_{x_2=0}^{\frac{2b(x_4-x_1)}{1-5b}}
\int_{x_3=0}^{\frac{2b(x_6-x_1)}{1-5b}}
\int_{x_5=x_4}^\infty
\int_{x_7=x_6}^\infty
I
dx_1 dx_2 dx_3 dx_4 dx_5 dx_6 dx_7
\]
where the integrand is
\[
I = b(x_4+x_5+x_6+x_7)e^{-x_1-x_2-x_3-x_4-x_5-x_6-x_7}
\]
which evaluates to
\[
\frac{1}{100}\frac{(73b-17)b^3}{(4b-1)^3}.
\]
The region on which $x_6>x_4$ gives
\[
\frac{1}{100}\frac{(73b-17)b^3}{(4b-1)(16b^2-8b+1)}.
\]
 The other
three parts can be expressed and evaluated similarly, giving a final result
of
\[
\ETb = \frac{726b^4-601b^3+245b^2-55b+5}{5(1-b)(4b-1)^2}.
\]
The result follows by simple computations.
\qed

Our final computation is for a nine-edge tree. Its lengthy proof uses the same
method that we used in the previous computations, so we omit it.
\begin{prop}
\lab{pAdmis}
Let $T$ be the biased tree on four degree-3 vertices and six leaf vertices
with bias $b$ on edges incident to leaves and bias
$(1-6b)/4$ on
the other edges, where
$b$ is the unique zero of
\[
2372895b^6-3013200b^5+1501416b^4-389232b^3+56016b^2-4224b+128
\]
on $0<b<1/6$.
Then
\[
\ETb = \frac{-2(-128+2448b-17856b^2+60372b^3-88938b^4+37665b^5)}{9(32-600b+4212b^2-13122b^3+15309b^4)}
\]
which lies in $(0.8696, 0.8697)$ and hence $c^*=0.8697$ is
$k$-admissible for $k=9$.
\end{prop}
Computer simulations suggest that the value of
$c^*$ decreases only slightly as $k$ is increased further  so we do not
pursue this here.

\section{Circumference of a random prekernel with given degree sequence}
In the previous sections we have established Lemma~\ref{lemmaHeavy},
an a.a.s.\ upper bound on the weight of the heaviest
cycle in certain randomly-edge-weighted pseudographs. In this section
we use that lemma
to establish an upper bound on the circumference of
a random prekernel whose degree sequence satisfies certain conditions.
In later sections we will see that the degree sequence of the prekernel
of $\Gnm$ a.a.s.\ satisfies these conditions,
 allowing us to use this result to establish an a.a.s.\
upper bound on the circumference of the prekernel of $\Gnm$.

One of the challenges in this section arises because Lemma~\ref{lemmaHeavy}
is a statement
about non-random pseudographs with random edge weightings, while we are
proving a statement about random prekernels. The kernel configuration
model of Pittel and Wormald, described below, allows us to rigorously
make this transition. It combines a pairing model, for generating the
kernel, with a random sequence of weights on the kernel edges.

Another challenge in this section is to show
that the conditions on the degree sequence imply
that the hypotheses of Lemma~\ref{lemmaHeavy} are satisfied.
One hypothesis requires that
there are few edges incident with vertices whose degree exceeds 3.
Another hypothesis requires that the number of short cycles in the kernel be small.
In {\L}uczak's proof of his upper bound for the circumference of
$\Gnm$ in the supercritical phase, he established the first hypothesis
by direct enumeration over degree sequences.
(See the proof of Theorem 10 in~\cite{Lcyc}.)
However, {\L}uczak does not require the second hypothesis, so we
will need to prove it here. We will see that, without much extra effort,
our proof of the
second hypothesis gives an alternative derivation of the first hypothesis.
In~\cite{boll} and \cite{W} there are results about short cycles arising
in this  pairing model. However,
 these results apply only when the maximum degree is bounded,
so they cannot be used for our application.

We are interested in studying prekernels with a given degree sequence
$\vd=(d_i)$. We say that $\vd$ is a \emph{prekernel degree sequence} if
its number of terms $v=v(\vd)$ is finite, each term is a positive integer
at least 2, and $r=r(\vd)=\sum_i(d_i-2)$ is even.
For $j = 2, 3, \ldots$ we define
\bel{eDj}
D_j=D_j(\vd)=|\{ i : d_i=j \}|.
\ee

The \emph{kernel configuration model} $\Hc(\vd)$ is used to generate prekernels
with degree sequence $\vd$.
It has been used successfully to calculate improved
estimates for the size of the core, excess, and tree mantle \cite{PW}.
We describe the model next.

For each $i$ with $d_i \geq 3$ create a set $S_i$ of $d_i$ points. Let $\P$
be the set of perfect matchings on the union of these sets of points and choose $P\in\P$
uniformly at random. Then, assign
the remaining numbers $\{ i : d_i=2 \}$ to the edges of the perfect
matching and, for each edge, choose a linear order for these numbers.
The assignments and the linear ordering, denoted by $f$, are chosen
uniformly at random. The pair $(P,f)$ defines a random configuration
in the model $\Hc(\vd)$.

Each configuration $(P,f)$ corresponds to a prekernel $G(P,f)$ by collapsing each set $S_i$ to a vertex
(producing a kernel $K(P)$) and placing the degree-2 vertices
on the edges of the kernel according to the assignment and linear orderings.

\begin{lemma}
\lab{lCircumDeg}
Let $\vd=\vd(n)$ be a prekernel degree sequence satisfying $v=v(\vd)\to\infty$,
$r=r(\vd)\to\infty$, $r=o(v)$, $D_3=D_3(\vd)\sim r$, and
\[
\sum_{i:d_i\ge 3} \binom{d_i}{2} < 4 r.
\]
Fix a positive integer $k\ge 2$ and suppose that the positive constant $c^*$ is $k$-admissible.
For a random configuration $(P,f)$ in $\Hc(\vd)$,
the longest
cycle in $G(P,f)$ has length a.a.s.\ at most
$c^* v$ as $n\to\infty$.
\end{lemma}
\proof
Define $\P^*$ to be the set of $P\in \P$ for which $K(P)$ has at most
$\sqrt{r}$ edges in cycles of length at most $k$. We will show that
a random configuration $(P,f)$ a.a.s.\ has $P\in\P^*$. Recall that $P$ is a random
perfect matching on the points in the union of the $S_i$. For $j\in \{1, 2, \ldots, k\}$,
the number of ways of choosing $j$ pairs of points to form a cycle
is at most
\[
\frac{1}{2j}
\left(
\sum_{i : d_i\ge 3} 2 \binom{d_i}{2}
\right)^j
=
O\left(r^j\right).
\]
The probability that $j$
given pairs of points appear in the pairing $P$ is
asymptotic to $$(\sum_{i:d_i\ge 3}d_i)^{-j}$$ since $j$ is bounded.
Now
\begin{eqnarray*}
\sum_{i:d_i\ge 3}d_i &>& \sum_{i:d_i\ge 3}(d_i-2) \\
&=& \sum_{i}(d_i-2) \\
&=& r
\end{eqnarray*}
so the expected number of cycles of length $j$ is
$O(r^j r^{-j})=O(1)$. Since $k$ is fixed, the expected number of edges in such
cycles is also $O(1)$. By Markov's inequality, the number of edges
in cycles of length
$j$ is a.a.s.\ bounded above by any function $\omega=\omega(n)\to\infty$,
in particular $\sqrt{r}/k$. Thus, a.a.s.\ $P\in\P^*$.

Let $(P,f)$ be a random configuration from $\Hc(\vd)$.
Define $G'(P,f)$ to be the edge-weighted pseudograph whose underlying
pseudograph is $K(P)$ and whose edge-weight on $e$, for each edge $e$,
is one more than the number of vertices assigned to $e$ by $f$.
Let $A$ be the event that the heaviest cycle in $G'(P,f)$ has weight at most
$c^* v$. Let $P_0$ be the $P^*$ minimizing
$\pr[A \mid P=P^*]$ over $P^*\in\P^*$. The minimum exists because $\P^*$ is
finite. Next we verify that, conditioned
on $P=P_0$, $G'(P,f)$ satisfies the hypotheses of
Lemma~\ref{lemmaHeavy}. The number of vertices $v'$ of $G'(P,f)$ is at least
$D_3\sim r\to \infty$. The minimum degree is at least 3 because it is a
kernel. The number of edges incident to cycles of length at most $k$
(including loops and parallel edges) is at most $\sqrt{r}=o(r)=o(v')$ since
$P_0\in\P^*$. The number of edges incident to vertices of degree greater than
3 is at most
\begin{eqnarray*}
\sum_{j:d_j\ge 4}d_j
&\le& 2\sum_{j:d_j\ge 4}(d_j-2) \\
&=&   2\sum_j (d_j-2) - 2D_3 \\
&=&   2r - 2D_3 \\
&=&   o(r)
\end{eqnarray*}
which is $o(v')$. The number of edges $m'$ satisfies
\begin{eqnarray*}
2m' &=&
\sum_{j:d_j\ge 3} d_j\\
&=& 3D_3 + o(r)
\end{eqnarray*}
by the previous calculation, so $m'=O(D_3)=O(r)=o(v)$. To see that
$m'$ is little-oh of the sum $N$ of the edge-weights,  observe that
$N$ is the number of edges of $G(P,f)$, which is $\sum_j d_j/2 \ge v$.
Next observe that the edge weights form a sequence of positive integers
that is determined by the assignment $f$ in the random configuration.
There are exactly $|\{i:d_i=2\}|!$ choices for $f$ that produce any given
sequence, so the sequence is chosen uniformly at random. We have shown that
the hypotheses of Lemma~\ref{lemmaHeavy} hold for $G'(P,f)$ conditioned
on $P=P_0$, so we have
$\pr[A \mid P=P_0]=1-o(1)$. Now
\begin{eqnarray*}
\pr[A] &\ge& \sum_{P^*\in\P^*} \pr[A \mid P=P^*] \pr[ P = P^* ] \\
&\ge& \pr[A \mid P=P_0]\sum_{P^*\in\P^*} \pr[ P = P^* ]
\end{eqnarray*}
by the choice of $P_0$. Since we showed $P\in\P^*$ a.a.s.\ we get
$\pr[A]=1-o(1)$; that is, the heaviest cycle in $G'(P,f)$ a.a.s.\
has weight at most $c^* v$. But if $C$ is a cycle in $G(P,f)$
of some length $l$, $C$ corresponds naturally to a cycle in $G'(P,f)$ of
weight $l$. So the longest cycle in $G(P,f)$ a.a.s.\ has length at most
$c^* v$.
\qed

\begin{cor}
\lab{corCircumDeg}
Let $\vd=\vd(n)$ be a prekernel degree sequence satisfying $v=v(\vd)\to\infty$,
$r=r(\vd)\to\infty$, $r=o(v)$, $D_3=D_3(\vd)\sim r$, and
\[
\sum_{i:d_i\ge 3} \binom{d_i}{2} < 4 r.
\]
Fix $k\ge 2$ and suppose that the positive constant $c^*$ is $k$-admissible.
Let $G$ be chosen uniformly at random from all prekernels with degree sequence
$\vd$.  The longest cycle in $G$ has length a.a.s.\ at most
$c^* v$ as $n\to\infty$.
\end{cor}
\proof
The probability space $\Hc(\vd)$, conditioned on the event that $G(P,f)$
is a simple graph, is a uniform probability space on the prekernels with degree
sequence $\vd$ (\cite{PW}, Lemma~3). By Lemma~5 in~\cite{PW}, $G(P,f)$ is
a.a.s.\ a simple graph. (In fact, Lemma~5 in~\cite{PW} is stated with
an additional
hypothesis on $\max d_i$, but this hypothesis is not used in the proof.)
The result now follows from Lemma~\ref{lCircumDeg}.
\qed

\section{Truncated multinomial distribution}
In order to apply Corollary~\ref{corCircumDeg} to the prekernel of
$\Gnm$, we must verify the hypotheses about properties of the degree sequence.
We give a new derivation of these properties, which will require some facts
about the following distribution.

Let $v$ and $t$ be positive integers. The probability space
$\multi(v,t)$ consists of vectors $(d_1, d_2, \ldots, d_v)$
with distribution
\[
\pr[ d_1=j_1, d_2=j_2, \ldots, d_v=j_v] =
\frac{t!}{v^t j_1! j_2! \cdots j_v!}
\]
for any vector $(j_1, j_2, \ldots, j_v)$ of nonnegative integers summing
to $t$. This is the well-known multinomial distribution, modelling the number
of balls in each bin when each of $t$ balls is tossed into one of $v$ bins,
independently and uniformly at random.
The space $\multi(v,t)|_{\ge 2}$ is obtained from $\multi(v,t)$ by
conditioning on the event that each $d_i\ge 2$.

\begin{lemma}
\lab{lemmaMulti}
Let $v=v(n)$ and $r=r(n)$ satisfy $v\to\infty$, $r\to\infty$ and
$r=o(v)$. If the random vector
$\vd$ is distributed as $\multi(v,2v+r)|_{\ge 2}$ then a.a.s.\
$D_3(\vd)\sim r$ and
\[
\sum_{i:d_i\ge 3} \binom{d_i}{2} < 4 r.
\]
\end{lemma}
\proof
Define the positive number $\lambda$ by
\[
\frac{ \lambda (e^\lambda - 1) }{ e^\lambda - 1 - \lambda } = 2 + \frac{r}{v}.
\]
In \cite{CW}, the authors show that $\lambda$ exists and they use it to define
a vector of independent truncated Poisson random variables which approximate
$\multi(v,2v+r)|_{\ge 2}$ as follows. Define the random variable $Y$ taking
values $j=2, 3, \ldots$ according to the distribution
\[
\pr[Y=j] = p_j = \frac{\lambda^j}{j!( e^\lambda - 1 - \lambda )}.
\]
Consider the probability space formed by vectors $\Yv=(Y_1, Y_2, \ldots, Y_v)$
of $v$ independent copies of
$Y$ and let $\Sigma$ be the event that their sum satisfies $\sum_i Y_i=2v+r$.
For nonnegative integers $j_1, j_2, \ldots, j_v$ summing to $2v+r$ with
each $j_i\ge 2$ we have
\[
\pr[Y_1=j_1, \ldots, Y_v=j_v] =
\frac{\lambda^{2v+r}}{(e^\lambda - 1 - \lambda)^v}
\prod_{i=1}^v \frac{1}{j_i!}
\]
so this probability space, conditioned on $\Sigma$, is identical to
$\multi(v,2v+r)|_{\ge 2}$.
Equation~(5.7) in~\cite{PW} says that
\[
\pr \left[
\max_i Y_i < \log v \mbox{ or }
\sum_{i:Y_i\ge 3} \binom{Y_i}{2} \ge 4 r \mid \Sigma
\right]
=
O(r^{-1} + rv^{-1}).
\]
The second claim in the lemma follows.
Theorem~4(a) in~\cite{PWae} states that for $r\to\infty$,
\[
\pr[\Sigma] = \frac{1+O(r^{-1})}{\sqrt{2\pi v c(1+\bar\eta-c)}}
\]
where $c(1+\bar\eta-c)\sim c-k = r/v$ by Equation~(20) in~\cite{PWae}.
It follows that
\bel{eSigma}
\pr[\Sigma]^{-1}=O(\sqrt{r}).
\ee
To establish the first claim in the lemma, we observe that $D_3(\Yv)$ is
distributed as a binomial random variable with $v$ trials and $p_3$ probability
of success.
By Chernoff's bound,
\[
\pr \left[ \left| D_3(\Yv) - vp_3 \right| > a \right]
<
2\exp( -a^2 / (3vp_3 ) )
\]
for $0 < a \le vp_3$.
Recalling \eqn{eSigma}, in $\multi(v,2v+r)|_{\ge 2}$ we have
\[
\pr \left[ \left| D_3(\vd) - vp_3 \right| > a \right]
=
O(\sqrt{r})
\exp( -a^2 / (3vp_3 ) ).
\]
Setting $a=\sqrt{vp_3} \log r$ (which satisfies $a\le vp_3$, as we will see
shortly) we get $D_3(\vd)=vp_3+O(\sqrt{vp_3}\log r)$
with probability $1-O(\exp(-(\log r)^2/3))$. Now $\lambda\sim 3rv^{-1}$ by
Theorem~1(a) in \cite{PWae}, so
\begin{eqnarray*}
vp_3 &=& v \frac{\lambda^3}{3!( e^\lambda - 1 - \lambda )} \\
&=& v \frac{\lambda^3}{3!(\lambda^2/2+O(\lambda^3))}\\
&=& \frac{1}{3} v \lambda (1+O(\lambda))\\
&=& r(1+O(rv^{-1}))
\end{eqnarray*}
giving us $D_3(\vd)\sim r$ a.a.s.\ as required.
\qed

\section{Properties of vertex degrees in $\Gnm$}
Now we proceed to establish the properties of the degree sequence of the
prekernel of $\Gnm$ that are required to apply Corollary~\ref{corCircumDeg}.
Recall that we are assuming $M=M(n)=n/2+s$ for some $s=s(n)$
satisfying $s=o(n)$ and $n^{2/3}=o(s)$.
For this range of $M$, it is well-known that $\Gnm$ a.a.s.\ has a unique
component with maximum number of vertices~\cite{B84}, which we call the largest
component.

We begin by showing a.a.s.\ there are few vertices in the core that lie outside
the largest component. The next result is part of the proof of Theorem~4 of \cite{Lcyc}. Here we
present a slightly more thorough proof.
\begin{lemma}
\lab{lCycCmp}
Let $M=M(n)=n/2+s$ for some $s=s(n)$ satisfying $s=o(n)$ and $n^{2/3}=o(s)$.
The number of vertices in cycles of $\Gnm$ not in the largest
component is a.a.s.\
at most $\omega n/s$ for any $\omega=\omega(n)\to\infty$.
\end{lemma}
\proof
Let $\bar G$ be the graph formed from $\Gnm$ by removing its (lexicographically
first) largest component.
Let $n(\bar G)$ and $M(\bar G)$ represent its number of vertices and edges,
respectively.
Let $\eps>0$ and define $S$ to be the set of ordered pairs $(\bar n, \bar M)$
satisfying
\begin{enumerate}
\item $(1-\eps)4s \le n - \bar n \le (1+\eps)4s$,
\item $(1-\eps)4s \le M - \bar M \le (1+\eps)4s$, and
\item $\pr[ n(\bar G) = \bar n, M(\bar G) = \bar M ] > 0$.
\end{enumerate}
It is known that the largest component of $\Gnm$ has a.a.s.\
$4s(1+o(1))$ vertices and $4s(1+o(1))$ edges~\cite{B84,L90}.
So a.a.s.\ $(n(\bar G), M(\bar G))\in S$. For $(\bar n, \bar M)\in S$ we have
$\bar M \le M - 4s(1-\eps) = n/2 + s - 4s(1-\eps)$ and
$n \le \bar n + 4s(1+\eps)$, giving us
\bel{eqBarMUpper}
\bar M \le \bar n / 2 - s(1-6\eps).
\ee
To estimate the number $X$ of vertices in cycles in $\bar G$ we let
$(\bar n, \bar M)\in S$ and condition on the non-empty event
$n(\bar G) = \bar n, M(\bar G) = \bar M$. In the conditioned space, $\bar G$
is equally likely to be any graph on $\bar n$ vertices and $\bar M$ edges.
For $3\le k \le \bar n$ the number of such graphs having a cycle of length
$k$ is at most
\[
\binom{\bar n}{k}\frac{k!}{2k}\binom{\binom{\bar n}{2}}{\bar M - k}
\]
so the expected value of $X$ in this conditioned space is
\begin{eqnarray*}
\ex[X\mid n(\bar G) = \bar n, M(\bar G) = \bar M]
&\le&
\sum_{k=3}^{\bar n} k \binom{\bar n}{k}\frac{k!}{2k}\binom{\binom{\bar n}{2}}{\bar M - k}\binom{\binom{\bar n}{2}}{\bar M}^{-1} \\
&=&
\frac{1}{2}\sum_{k=3}^{\bar n}\frac{\bar n !}{(\bar n - k)!}
                              \frac{\left(\binom{\bar n}{2}-M-k\right)!}{\left(\binom{\bar n}{2}-M\right)!}
                              \frac{\bar M !}{(\bar M - k)!}\\
&<&
\frac{1}{2}\sum_{k=3}^{\bar n} \frac{\bar{n}^k \bar{M}^k}{(\binom{\bar n}{2}-\bar n)^k}.
\end{eqnarray*}
Using \eqn{eqBarMUpper} this becomes
\begin{eqnarray*}
\ex[X\mid n(\bar G) = \bar n, M(\bar G) = \bar M]
&<&
\frac{1}{2}\sum_{k=3}^{\bar n} \left( \frac{\frac{\bar n}{2}-s(1-6\eps)}{\frac{\bar n-1}{2}-1}\right)^k \\
&<&
\frac{1}{2}\sum_{k=3}^{\bar n}\left( 1 - \frac{s(1-6\eps)}{\bar n - 3}\right)^k\\
&<&
\frac{1}{2}\sum_{k=3}^{\infty}\left( 1 - \frac{s(1-6\eps)}{\bar n - 3}\right)^k\\
&=&
\frac{1}{2}\times\frac{\bar n - 3}{s(1-6\eps)}
\end{eqnarray*}
which is at most $n/s$ for $n$ sufficiently large.
By Markov's inequality,
\begin{eqnarray*}
\pr[X \ge \omega n/s \mid n(\bar G) = \bar n, M(\bar G) = \bar M]
&\le&
\frac{\ex[X\mid n(\bar G) = \bar n, M(\bar G) = \bar M]}{\omega n/s}\\
&<&\frac{1}{\omega}.
\end{eqnarray*}
So
\begin{eqnarray*}
\pr[X < \omega n/s]
&\ge&
\pr[X < \omega n/s, (n(\bar G), M(\bar G))\in S] \\
&=&
\sum_{(\bar n, \bar M)\in S}
\big(
\pr[ X < \omega n/s \mid n(\bar G) = \bar n, M(\bar G) = \bar M]\\
&&\times
\pr[ n(\bar G) = \bar n, M(\bar G) = \bar M]\big)\\
&\ge&
\left(1-\frac{1}{\omega}\right)
\pr[(n(\bar G), M(\bar G))\in S]  \\
&=&
\left(1-\frac{1}{\omega}\right)
(1-o(1))
\end{eqnarray*}
since a.a.s.\
$(n(\bar G), M(\bar G))\in S$.
Therefore a.a.s.\ $X < \omega n/s$.
\qed

Instead of proving results about the degree sequence of the prekernel of
$\Gnm$ directly, we will actually prove results about the degree sequence
of the core. The next result will allow us to transfer results about the core to
the prekernel.
\begin{lemma}
\lab{lConstCore}
Let $M=M(n)=n/2+s$ for some $s=s(n)$ satisfying $s=o(n)$ and $n^{2/3}=o(s)$.
The core of the largest component of $\Gnm$
is a.a.s.\ formed from the core of $\Gnm$ by removing $o(s^2/n)$ vertices of degree 2.
Also, the prekernel
of $\Gnm$
is a.a.s.\ formed from the core of $\Gnm$ by removing $o(s^2/n)$ vertices of
degree 2.
\end{lemma}
\proof
It is well-known that the largest component of $\Gnm$ is a.a.s.\ the only
component that has more than one cycle. (See Theorem~5.12 in~\cite{JLR}.)
So, the core of $\Gnm$ is a.a.s.\ composed of the core of the largest component
together with some cycle components. By Lemma~\ref{lCycCmp}
the number of vertices in the cycle components is a.a.s.\ at most
$\omega n/s = (s/n^{2/3})(n/s) = n^{1/3} = (n^{2/3})^2/n = o(s^2/n)$,
since we may take $\omega=s/n^{2/3}$. Because the largest component of $\Gnm$
a.a.s.\ contains more than one cycle, it follows that these cycle components
are a.a.s.\ all of the cycle components in the core of $\Gnm$, making the
prekernel a.a.s.\ equal to the core of the largest component.
\qed

Now we establish the required properties of the degree sequence of
the prekernel of $\Gnm$.
Recall the definitions of $D_j(\vd)$, $v(\vd)$, $r(\vd)$ from \eqn{eDj}.
The results about $r(\vd)$, $v(\vd)$, and $D_3(\vd)$ in the following lemma
were used by {\L}uczak~\cite{Lcyc}. He used the estimate of
$r(\vd)$ from~\cite{L90} to establish
estimates for $D_j(\vd)$ and $v(\vd)$ by direct enumeration over degree
sequences. Instead of studying the prekernel directly, he studied the
core of the largest component.
Our proof method is different, using the known estimates of
$v(\vd)$ and $r(\vd)$ to establish results about the degree sequence of
the core.
\begin{lemma}
\lab{lPrekDegSeq}
Let $M=M(n)=n/2+s$ for some $s=s(n)$ satisfying $s=o(n)$ and $n^{2/3}=o(s)$.
Let $\vd$ be the degree sequence of the prekernel of $\Gnm$. Then, a.a.s.\
$v(\vd)\sim 8s^2/n$, $r(\vd)\sim D_3(\vd)\sim 32s^3/(3n^2)$, and
\[
\sum_{i : d_i\ge 3} 2 \binom{d_i}{2} < 4r(\vd).
\]
\end{lemma}
\proof
By Lemma~\ref{lConstCore} the prekernel differs from the
core a.a.s.\ by $o(s^2/n)$ vertices of degree 2. Thus, it suffices to prove
the lemma for the degree sequence $\vd$ of the core. Appealing to
Lemma~\ref{lConstCore} again, the core differs from the core of the largest
component a.a.s.\ by $o(s^2/n)$ vertices of degree 2. It is known~\cite{PW}
 that the degree sequence $\vd'$ of the core of the largest component a.a.s.\
has $v(\vd')\sim 8s^2/n$ and $r(\vd')\sim 32s^3/(3n^2)$, so we must have
a.a.s.\ $v(\vd)\sim 8s^2/n$ and $r(\vd)\sim 32s^3/(3n^2)$ also.
Letting $\eps>0$, this means $\vd \in S$ a.a.s.\ where $S$ is the set of
ordered pairs $(\bar v,\bar r)$ satisfying
\begin{enumerate}
\item $(1-\eps)8s^2/n       \le  \bar v  \le (1+\eps)8s^2/n$,
\item $(1-\eps)32s^3/(3n^2) \le  \bar r  \le (1+\eps)32s^3/(3n^2)$, and
\item $\pr[ v(\vd)=\bar v, r(\vd)=\bar r ] > 0$.
\end{enumerate}
We note that for $(\bar v, \bar r)\in S$ we have
$\bar r = o(\bar v)$ since $(s^3/n^2)/(s^2/n)=s/n=o(1)$
and both $\bar v\to \infty$ and $\bar r\to \infty$ since $n^{2/3}=o(s)$.

To establish the remaining properties of the degree sequence of the core of
$\Gnm$ we use Theorem~2 of~\cite{CW},
which proves the existence of a probability space of ordered pairs $(G,I)$
in which
\begin{enumerate}
\item $G$, conditioned on the event $I=1$, is distributed as the core of $\Gnm$,
\item $\pr[I=1] = \Omega(1)$, and
\item the degree sequence $\vd(G)$ of $G$, conditioned on $v(\vd(G))=\bar v$
and $r(\vd(G))=\bar r$, is distributed as
$\multi(\bar v, 2\bar{v}+\bar{r})|_{\ge 2}$.
\end{enumerate}
(The statement of Theorem~2 of~\cite{CW} actually includes the hypothesis
$M\ge n$ which is not satisfied here; however, that hypothesis is not
needed for their proof.)

Write $v=v(\vd(G))$, $r=r(\vd(G))$ and let $A$ be the event that
\begin{enumerate}
\item $(1-\eps)32s^3/(3n^2) \le D_3(\vd(G)) \le (1+\eps)32s^3/(3n^2)$, and
\item $\sum_{i : d_i(G) \ge 3} 2 \binom{d_{i}(G)}{2} < 4r$.
\end{enumerate}
To prove the lemma, we must show $\pr[A\mid I=1]=1+o(1)$ or equivalently
$\pr[A^C \mid I=1]=o(1)$, where $X^C$ denotes the complement of event $X$.
We begin by writing
\[
\pr[A^C \mid I=1] = \pr[A^C, (v,r)\in S \mid I=1]
                  + \pr[A^C, (v,r)\not\in S \mid I=1].
\]
The second term is at most $\pr[(v,r)\not\in S \mid I=1]$
which is $o(1)$ because $\vd(G)$, conditioned on $I=1$, is distributed like
the degree sequence of the core of $\Gnm$. We write the first term as
\begin{eqnarray*}
\pr[A^C, (v,r)\in S \mid I=1]
&=&
\sum_{(\bar v,\bar r)\in S} \pr[A^C, v=\bar v, r=\bar r \mid I=1] \\
&=&
\pr[I=1]^{-1} \sum_{(\bar v,\bar r)\in S} \pr[A^C, v=\bar v, r=\bar r, I=1] \\
&\le&
\pr[I=1]^{-1} \sum_{(\bar v,\bar r)\in S} \pr[A^C, v=\bar v, r=\bar r] \\
&=&
\pr[I=1]^{-1} \sum_{(\bar v,\bar r)\in S} \frac{\pr[A^C, v=\bar v, r=\bar r] \pr[v=\bar v, r=\bar r]}{\pr[v=\bar v, r=\bar r]} \\
&=&
\pr[I=1]^{-1} \sum_{(\bar v,\bar r)\in S} \pr[A^C \mid v=\bar v, r=\bar r] \pr[v=\bar v, r=\bar r]\\
&\le&
\pr[I=1]^{-1} \pr[A^C \mid v=\hat v, r=\hat r] \pr[(v, r)\in S]
\end{eqnarray*}
where $(\hat v, \hat r)$ is the ordered pair maximizing
$\pr[A^C \mid v=\bar v, r=\bar r]$ over all
$(\bar v, \bar r)\in S$. (The maximum exists because $S$ is finite.)
We now use the properties of the distribution of $(G,I)$ to  estimate
each of $\pr[I=1]^{-1}$, $\pr[A^C \mid v=\hat v, r=\hat r]$, and $\pr[(v, r)\in S]$.
We have already noted that $\pr[I=1]=\Omega(1)$, so we have
$\pr[I=1]^{-1}=O(1)$.
Since $(\hat v, \hat r)\in S$ we have $\hat v \to \infty$, $\hat r \to \infty$, and $\hat r = o(\hat v)$. We know that, conditioned on
$v=\hat v$ and $r=\hat r$, the degree sequence of $G$ is
 distributed as $\multi(\hat v, 2\hat{v}+\hat{r})|_{\ge 2}$.
Lemma~\ref{lemmaMulti} tells us that event $A$ occurs a.a.s.\  in this model
 so we have $\pr[A^C \mid v=\hat v, r=\hat r]=o(1)$. Finally,
we may crudely estimate $\pr[(v, r)\in S]=O(1)$.
Combining these estimates we get $\pr[A^C \mid I=1]=o(1)$ as required.
\qed

\section{Circumference of $\Gnm$}
\begin{lemma}
\lab{lemmaCircumGnm}
Let $M=M(n)=n/2+s$ for some $s=s(n)$ satisfying $s=o(n)$ and $n^{2/3}=o(s)$.
Fix $k$ and suppose that the positive constant $c^*$ is $k$-admissible.
The circumference of $\Gnm$ is a.a.s.\ at most $(8 c^*+o(1)) s^2/n$.
\end{lemma}
\proof
Every cycle in a graph lies in the graph's core.
By Lemma~\ref{lConstCore} the prekernel $G$ of $\Gnm$ is formed from the core
of $\Gnm$ by removing $o(s^2/n)$ vertices of degree 2. So, to prove the lemma,
it suffices to show that the circumference of $G$
is a.a.s.\ at most $(8 c^* + o(1)) s^2/n$.

By Lemma~\ref{lPrekDegSeq} there exists $\omega=\omega(n)\to\infty$ such
that the degree sequence $\vd(G)$ of $G$ a.a.s.\ lies in
the set $D$ of prekernel degree sequences $\vd$ satisfying
\begin{enumerate}
\item
\[
\sum_{i : d_i\ge 3} 2 \binom{d_i}{2} < 4r(\vd),
\]
\item
$(1-\omega^{-1})8s^2/n \le v(\vd) \le (1+\omega^{-1}) 8s^2/n$,
\item
$(1-\omega^{-1})32s^3/(3n^2) \le r(\vd) \le (1+\omega^{-1}) 32s^3/(3n^2)$,
\item
$(1-\omega^{-1})32s^3/(3n^2) \le D_3(\vd) \le (1+\omega^{-1}) 32s^3/(3n^2)$,
and
\item
$\pr[ \vd(G) = \vd ] > 0$.
\end{enumerate}
Define $A$ to be the event that the circumference of $G$ is at most
$c^* v(\vd(G))$. We have
\[
\pr[A] \ge \sum_{\vd\in D} \pr[A\mid \vd(G)=\vd] \pr[\vd(G)=\vd].
\]
Suppose that $\pr[A\mid \vd(G)=\vd]$ is minimized over $\vd\in D$ by
$\vd=\vd^*$. (The minimum exists since $D$ is finite.) Then
\begin{eqnarray*}
\pr[A] &\ge& \pr[A\mid \vd(G)=\vd^*] \sum_{\vd\in D}  \pr[\vd(G)=\vd]\\
&=& \pr[A\mid \vd(G)=\vd^*] (1+o(1))
\end{eqnarray*}
since a.a.s.\ $\vd\in D$.

In general, the number of graphs on $n$ vertices and $M$ edges that have a
given graph as their prekernel depends only on the number of vertices and
edges of the given prekernel. So, conditioning on the event $\vd(G)=\vd^*$,
$G$ is equally likely to be each prekernel with degree sequence $\vd^*$.
The probability $P[A\mid \vd(G)=\vd^*]$ is thus the probability that
a graph, chosen uniformly at random from all prekernels of degree sequence
$\vd^*$,  has circumference a.a.s.\ at most $c^* v$. Since $\vd^*\in D$,
we may apply Corollary~\ref{corCircumDeg} to conclude that this
probability is $1+o(1)$. Thus $\pr[A]=1+o(1)$. In other words, the circumference
of $G$ is a.a.s.\ at most $c^* v(\vd(G))$. But we have seen
$v(\vd(G))\sim 8s^2/n$ a.a.s.\ so the circumference is a.a.s.\ at most
$(8c^* + o(1))s^2/n$, as required.
\qed

{\bf Proof of Theorem~\ref{tMain}.}
Proposition~\ref{pAdmis} tells us that $c^*=0.8697$ is $k$-admissible for $k=9$.
By Lemma~\ref{lemmaCircumGnm} the circumference of $\Gnm$ is a.a.s.\ at most
$(8 c^*+o(1)) s^2/n < (6.958+o(1)) s^2/n$.
\qed

\newcommand\AAP{\emph{Adv. Appl. Probab.} } 
\newcommand\JAP{\emph{J. Appl. Probab.} }
\newcommand\JAMS{\emph{J. \AMS} }
\newcommand\MAMS{\emph{Memoirs \AMS} }
\newcommand\PAMS{\emph{Proc. \AMS} }
\newcommand\TAMS{\emph{Trans. \AMS} }
\newcommand\AnnMS{\emph{Ann. Math. Statist.} }
\newcommand\AnnPr{\emph{Ann. Probab.} }
\newcommand\CPC{\emph{Combin. Probab. Comput.} }
\newcommand\JMAA{\emph{J. Math. Anal. Appl.} }
\newcommand\RSA{\emph{Random Structures Algorithms} }
\newcommand\ZW{\emph{Z. Wahrsch. Verw. Gebiete} }
\newcommand\DMTCS{\jour{Discr. Math. Theor. Comput. Sci.} }

\newcommand\AMS{Amer. Math. Soc.}
\newcommand\Springer{Springer-Verlag}
\newcommand\Wiley{Wiley}

\newcommand\vol{\textbf}
\newcommand\jour{\emph}
\newcommand\book{\emph}
\newcommand\inbook{\emph}
\def\no#1#2,{\unskip#2, no. #1,} 

\newcommand\webcite[1]{\hfil\penalty0\texttt{\def~{\~{}}#1}\hfill\hfill}
\newcommand\webcitesvante{\webcite{http://www.math.uu.se/\~{}svante/papers/}}
\newcommand\arxiv[1]{arXiv:#1}

\def\nobibitem#1\par{}

\end{document}